\documentclass[a4paper]{IEEEtran}
\usepackage[latin1]{inputenc}
\usepackage[cmex10]{amsmath}
\usepackage{amsfonts}
\usepackage{amssymb}
\usepackage{graphicx}

\usepackage{verbatim}
\usepackage{array}
\usepackage{multirow}
\usepackage{dcolumn}
\usepackage{color}
\usepackage[noadjust]{cite}
\usepackage{url}
\usepackage{balance}
\usepackage[usenames,dvipsnames]{xcolor}
\usepackage{booktabs}
\usepackage{amsmath}
\usepackage{nomencl}
\usepackage{breqn}
\makenomenclature

\usepackage{etoolbox}
\renewcommand\nomgroup[1]{%
  \item[\bfseries
  \ifstrequal{#1}{S}{Sets}{%
  \ifstrequal{#1}{I}{Index}{%
  \ifstrequal{#1}{C}{Constants}{%
  \ifstrequal{#1}{V}{Decision Variables}{}}}}%
]}

\DeclareGraphicsExtensions{.eps}

\begin{document}
\title{Optimal Generation and Transmission Expansion Planning Addressing Short-Term Constraints with Co-optimization of Energy and Reserves}

\author{Alessandro~Soares,~Ricardo~Perez,~Fernanda~Thome

\thanks{Alessandro Soares, Ricardo Perez, Fernanda Thome are from PSR, Rio de Janeiro, RJ, Brazil (e-mail: \mbox{alessandro@psr-inc.com}; \mbox{ricardo@psr-inc.com}; \mbox{fernanda@psr-inc.com})}
}

\maketitle


\begin{abstract}
The penetration of variable renewable energy (VRE) in electrical systems has changed the way the expansion planning is treated. This kind of resource has great variability in small amounts of time, which makes it important to represent hourly constraints that requires chronology. Besides that, the generation reserves should also be adjusted in order to capture the intermittent effect, and since many countries use rapid thermal plants as part of these reserves, unit commitment and ramp constraint have also become more significant. In this paper we propose a MILP expansion planning model that can represent hourly time steps while maintaining reasonable computational times, where both investment and operation problems are simultaneously solved. Because the planning horizons are long (decades), the resolution of the entire horizon in a single optimization problem would be computationally infeasible for large real systems, making it necessary, therefore, to apply a horizon decomposition heuristic in smaller sub-horizons, and use the representation of typical days and seasons to reduces the size of the problem.
\end{abstract}

\section*{Nomenclature}
\label{sec:nomenclature}

\subsection*{Constants}

\begin{IEEEdescription}[\IEEEusemathlabelsep\IEEEsetlabelwidth{$MMM$}]
\item[$w_\zeta^k$]{ Installed capacity, firm energy or firm capacity of generic project $\zeta$ relative to constraint $k$.}
\item[$\underline{w_k}$]{ Lower-bound of min/max constraint $k$.}
\item[$\overline{w_k}$]{ Upper-bound of min/max constraint $k$.}
\item[$w_\zeta^k$]{ Installed capacity, firm energy or firm capacity of generic project $\zeta$ relative to min/max constraint $k$.}
\item[$\underline{w_k}$]{ Minimum value (RHS) of min/max constraint $k$.}
\item[$\overline{w_k}$]{ Maximum value (RHS) of min/max constraint $k$.}
\item[$\overline{g_j}$]{ Maximum generation of thermal plant $j$}
\item[$\underline{g_j}$]{ Minimum generation of thermal plant $j$}
\item[$\rho_i$]{ Mean production factor of hydro plant $i$}
\item[$\overline{V_i}$]{ Maximum storage of hydro plant $i$}
\item[$\underline{V_i}$]{ Minimum storage of hydro plant $i$}
\item[$\overline{U_i}$]{ Maximum turbining of hydro plant $i$}
\item[$\overline{g_i}$]{ Maximum generation of hydro plant $i$}
\item[$\underline{u_i}$]{ Minimum turbining of hydro plant $i$}
\item[$\underline{q_i}$]{ Minimum total outflow of hydro plant $i$}
\item[$D_{t,d}$]{ Duration of typical day $d$ at season $t$}
\item[$\psi_{l,t,d,h,s}$]{ Renewable generation scenario for plant $l$, season $t$, typical day $d$, hour of the day $h$ and scenario $s$}
\item[$\eta^{+}_b$]{ Charge efficiency of battery $b$}
\item[$\eta^{-}_b$]{ Discharge efficiency of battery $b$}
\item[$\overline{V_b}$]{ Maximum storage of battery $b$}
\item[$\overline{q^{+}_b}$]{ Maximum charge capacity of battery $b$}
\item[$\overline{q^{-}_b}$]{ Maximum discharge capacity of battery $b$}
\item[$\overline{f_k^{+}}$]{ From$\rightarrow$To maximum flow capacity in transmission line $k$.}
\item[$\overline{f_k^{-}}$]{ To$\rightarrow$From maximum flow capacity in transmission line $k$.}
\item[$\gamma_k$]{ Susceptance of transmission line $k$.}
\item[$M$]{ Disjunctive constant.}
\item[$\overline{Imp_a}$]{ Maximum import amount of area $a$.}
\item[$\underline{Imp_a}$]{ Minimum import amount of area $a$.}
\item[$\overline{Exp_a}$]{ Maximum export amount of area $a$.}
\item[$\underline{Exp_a}$]{ Maximum import amount of area $a$.}
\item[$\underline{g_c}$]{Minimum value of generation constraint $c$}
\item[$\overline{g_c}$]{Maximum value of generation constraint $c$}
\item[$R_{c,t,d,h,s}$]{ Reserve requirement $c$, season $t$, typical day $d$, hour of the day $h$ and scenario $s$}
\item[$D_{n,t,d,h,s}$]{Inelastic demand associated to bus $n$, season $t$, typical day $d$, hour of the day $h$ and scenario $s$}
\item[$p_s$]{Probability of scenario $s$}
\item[$rt$]{Discount rate at the season}
\item[$co_j$]{The operation cost of thermal plant $j$}
\item[$sc_j$]{The start-up cost of thermal plant $j$}
\item[$co_i$]{O\&M cost of hydro plant $i$}
\item[$c_{\delta_i^v}$]{Minimum storage violation penalty of hydro plant $i$}
\item[$c_{\delta_i^u}$]{Minimum turbining violation penalty of hydro plant $i$}
\item[$c_{\delta_i^q}$]{Total outflow violation penalty of hydro plant $i$}
\item[$c_{\delta_c^G}$]{Violation penalty of generation constraint $c$}
\item[$c_{\delta_c^R}$]{Violation penalty of reserve requirement constraint $c$}
\item[$c_d$]{Deficit cost}
\item[$P_n^E$]{Elastic demand price of bus $n$}
\item[$ci_j$]{Investment cost of thermal project $j$}  
\item[$ci_i$]{Investment cost of hydro project $i$}  
\item[$ci_l$]{Investment cost of renewable project $l$}  
\item[$ci_b$]{Investment cost of battery project $b$}  
\item[$ci_k$]{Investment cost of transmission line project $k$} 

\end{IEEEdescription}

\subsection*{Sets}
\begin{IEEEdescription}[\IEEEusemathlabelsep\IEEEsetlabelwidth{$MMM$}]

\item[$R^{pre}$]{ Set of precedence constraints.}
\item[$P_k^{pre}$]{ Set of projects that belong to precedence constraint $k$}
\item[$R^{ctr}$]{ Set of minimum / maximum constraints.}
\item[$P_k^{ctr}$]{ Set of projects that belong to the constraint $k$.}
\item[$R^{ex}$]{ Set of exclusivity constraints.}
\item[$P_k^{ex}$]{ Set of projects that belong to exclusivity constraint $k$.}
\item[$R^{as}$]{ Set of association constraints.}
\item[$P_k^{as}$]{ Set of projects that belong to association constraint $k$.}
\item[$R^{ctr}$]{ Set of min/max constraints.}
\item[$P_k^{ctr}$]{ Set of projects that belong to the constraint $k$.}
\item[$M_{i}$]{ The set of plants upstream of hydro plant $i$}
\item[$K^p$]{ Set of circuits projects.}
\item[$K_a^{+}$]{ Set of transmission lines that arrive at area $a$ (To bus is in the area $a$ and the From bus is in a different area)}
\item[$K_a^{-}$]{ Set of transmission lines that leave at area $a$ (From bus is in the area $a$ and the To bus is in a different area)}
\item[$J_c^G$]{ Set of thermal plants that belongs to generation constraint $c$}
\item[$I_c^G$]{ Set of hydro plants that belongs to generation constraint $c$}
\item[$J_c^R$]{ Set of thermal plants that belongs to reserve constraint $c$}
\item[$I_c^R$]{ Set of hydro plants that belongs to reserve constraint $c$}
\item[$B_c^R$]{ Set of batteries that belongs to reserve constraint $c$}
\item[$J_n$]{Set of thermal plants that belong to bus $n$}
\item[$I_n$]{Set of hydro plants that belongs to bus $n$}
\item[$B_n$]{Set of batteries that belongs to bus $n$}
\item[$L_n$]{Set of renewable plants that belongs to bus $n$}
\item[$K_n^{+}$]{Set of transmission lines that arrive at bus $n$}
\item[$K_n^{-}$]{Set of transmission lines that leave bus $b$}
\item[$J_x$]{Set of thermal projects}  
\item[$I_x$]{Set of hydro projects}  
\item[$L_x$]{Set of renewable projects}  
\item[$B_x$]{Set of battery projects}  
\item[$K_x$]{Set of transmission line projects}  

\end{IEEEdescription}

\subsection*{Decision Variables}

\begin{IEEEdescription}[\IEEEusemathlabelsep\IEEEsetlabelwidth{$MMM$}]
\item[$x_\omega$]{ Decision variable of generic project $\omega$}
\item[$x_\zeta$]{ Decision variable of generic project $\zeta$}
\item[$\gamma_{j,t,d,h,s}$]{ Commitment decision of thermal plant $j$, season $s$, typical day $d$, hour of the day $h$ and scenario $s$}
\item[$g_{j,t,d,h,s}$]{ Generation decision of thermal plant $j$, season $s$, typical day $d$, hour of the day $h$ and scenario $s$}
\item[$x_j$]{Investment decision of thermal plant $j$}
\item[$st_{j,t,d,h,s}$]{Startup decision of thermal plant $j$, season $s$, typical day $d$, hour of the day $h$ and scenario $s$}
\item[$v_{i,t,s}$]{ Storage of hydro plant $i$, season $t$ and scenario $s$}
\item[$u_{i,t,s}$]{ Turbining of hydro plant $i$, season $t$ and scenario $s$}
\item[$s_{i,t,s}$]{ Spilling of hydro plant $i$, season $t$ and scenario $s$}
\item[$a_{i,t,s}$]{ Lateral streamflow arriving at hydro plant $i$, season $t$ and scenario $s$}
\item[$g_{i,t,d,h,s}$]{ Generation decision of hydro plant $i$, season $t$, typical day $d$, hour of the day $h$ and scenario $s$}
\item[$\delta^v_{i,t,s}$]{ Minimum storage violation decision of hydro plant $i$, season $t$ and scenario $s$}
\item[$\delta^u_{i,t,s}$]{ Minimum turbining violation decision of hydro plant $i$, season $t$ and scenario $s$}
\item[$\delta^1_{i,t,s}$]{ Minimum total outflow violation decision of hydro plant $i$, season $t$ and scenario $s$}
\item[$g_{l,t,d,h,s}$]{ Generation decision of renewable plant $l$, season $t$, typical day $d$, hour of the day $h$ and scenario $s$}
\item[$x_l$]{ Investment decision of renewable plant $l$}
\item[$v_{b,t,d,h,s}$]{ Storage of battery $b$, season $t$, typical day $d$, hour of the day $h$ and scenario $s$}
\item[$q^{+}_{b,t,d,h,s}$]{ Charge of battery $b$, season $t$, typical day $d$, hour of the day $h$ and scenario $s$}
\item[$q^{-}_{b,t,d,h,s}$]{ Discharge of battery $b$, season $t$, typical day $d$, hour of the day $h$ and scenario $s$}
\item[$x_b$]{ Investment decision of battery $b$}
\item[$f_{k,t,d,h,s}^{+}$]{From$\rightarrow$To flow in transmission line $k$, season $t$, typical day $d$, hour of the day $h$ and scenario $s$}
\item[$f_{k,t,d,h,s}^{-}$]{To$\rightarrow$From flow in transmission line $k$, season $t$, typical day $d$, hour of the day $h$ and scenario $s$}
\item[$x_k$]{Investment decision of transmission line $k$}
\item[$\theta_{b_k^{+},t,d,h,s}$]{ Nodal angle of bus $b_k^{+}$ (From bus of transmission line $k$.}
\item[$\theta_{b_k^{-},t,d,h,s}$]{ Nodal angle of bus $b_k^{-}$ (To bus of transmission line $k$.}
\item[$\delta^g_{c,t,d,h,s}$]{ Violation decision of generation constraint $c$, season $t$, typical day $d$, hour of the day $h$ and scenario $s$}
\item[$\delta^R_{c,t,d,h,s}$]{ Violation decision of reserve requirement constraint $c$, season $t$, typical day $d$, hour of the day $h$ and scenario $s$}
\item[$r_{j,t,d,h,s}$]{ Reserve allocated by thermal plant $j$, season $t$, typical day $d$, hour of the day $h$ and scenario $s$}
\item[$r_{i,t,d,h,s}$]{ Reserve allocated by hydro plant $i$, season $t$, typical day $d$, hour of the day $h$ and scenario $s$}
\item[$r_{b,t,d,h,s}$]{ Reserve allocated by battery $b$, season $t$, typical day $d$, hour of the day $h$ and scenario $s$}
\item[$DE_{n,t,d,h,s}$]{Elastic demand associated to bus $n$, season $t$, typical day $d$, hour of the day $h$ and scenario $s$}
\item[$\delta_{n,t,d,h,s}$]{Deficit at bus $n$, season $t$, typical day $d$, hour of the day $h$ and scenario $s$}

\end{IEEEdescription}

\section{Introduction}
\label{sec:introduction}

\IEEEPARstart{T}{he} interest on optimal power system expansion planning has increased worldwide. In developing countries of Latin America, Asia and Africa, with high load growth and limited financial resources, the emphasis is on the most cost-effective expansion plan \cite{ondraczek2014we,deichmann2011economics,kaygusuz2012energy}. In developed countries, load growth is usually flat. In these cases, Renewable Energy Sources (RES) are being built as part of decarbonization policies and to displace more expensive thermal plants \cite{haller2012decarbonization,oberthur2015decarbonization,capros2014european,ekins2004step}. In both cases, selecting the "best" of a group of alternatives is what characterizes the combinatorial nature of the expansion planning problem.

The main objective of the expansion planning process is to guarantee an appropriate balance between electricity supply and demand, i.e. to determine the optimal set of generating plants and transmission routes that should be constructed to meet the demand requirements along a study horizon (mid and long term), while minimizing a cost function considering: (i) investment (capital) and operation (fuel, O\&M, etc.) costs of generation plants and (ii) penalties of energy not supplied, also called deficit costs. 

In general terms, this decision process involves meeting economic, reliability and environmental criteria, within the framework of national policies on energy. In addition, one of the greatest challenges is how to deal with the uncertainties inherent in the planning process, such as the load growth, the hydrological inflows and the generation availability, especially in renewable based systems. Taken all the aforementioned facts into account, the expansion planning problem is modeled as a large and complex mixed integer multistage stochastic problem that must be solved by specialized optimization algorithm.

This paper presents a description of the methodology associated with the OptGen model \cite{optgenpsr}, a commercial computational tool for energy systems expansion planning, where two "Solution Strategies" are available: 
\begin{itemize}
    \item The benders decomposition strategy, proposed in \cite{campodonico2003expansion}: A decomposition of the investment and operation problem, where the master is a MILP investment problem and the slave is a multistage stochastic optimization of the operational problem that is solved using the SDDP algorithm, first proposed in \cite{pereira1991multi};
    \item The co-optimization strategy, which is the methodology described in this paper
\end{itemize}

The main characteristics of the model are:
\begin{itemize}
    \item Study horizons from 1 year up to several decades;
    \item Many different candidate projects may be contemplated in the study, such as:
    \begin{itemize}
        \item Production components: hydro, thermal and renewable plants (wind, solar, biomass, etc.);
        \item Interconnection links and transmission circuits (lines, transformers, DC links etc.);
        \item Storage devices such as Batteries, Hydro pump stations, etc.
        \item Gas pipelines, production nodes, regasification stations.
    \end{itemize}
    \item Detailed project's financial data, such as, investment costs, payment schedules, life-time, construction time;
    \item Detailed project specific data, such as, decision type (obligatory or optional), decision variable type (binary, integer or continuous), maximum and minimum en-trance dates, generating unit entrance schedule, etc.;
    \item Additional constraints, such as, firm energy/capacity constraints, exclusivity, association and precedence between projects, minimum and maximum additional capacity, generation capacity targets and so on;
    \item Unit commitment constraints
    \item Ramping constraints
    \item Co-optimization of energy and reserves
\end{itemize}

In summary, the objective of OptGen is to determine a least-cost investment schedule for the construction of new plant capacity (hydro, thermal and renewable projects), regional inter-connections (or detailed transmission circuits), gas production sources and gas pipelines. This is obtained by optimizing the trade-off between investment costs to build new projects and the expected value of operating and deficit costs.

This paper is organized as follows. In the Section \ref{sec:litreview}, a review of the current state-of-the-art policies, methodologies and models regarding systems with a high level of renewable penetration is presented. In Section \ref{sec:solmet}, we discuss the assumptions used by the methodology in order to make it computational tractable. In Section \ref{sec:uncertain}, we analyze how the uncertainties are taken into account in the proposed model. In Section \ref{sec:problemform}, we provide a detailed formulation of the proposed methodology. Finally, in Section \ref{sec:conclusion}, the final conclusions are presented.


\section{Literature review}\label{sec:litreview}

The increasing economic competitiveness of wind and solar generation sources, also called variable renewable energy sources (VRE), has widely studies in the literature. These energy resources reduces green-house gas emissions, as studied in \cite{renewableenergyeurope}. Besides that, \cite{irenageopolitical} showed that in a renewable energy economy, since renewable energy potential is available everywhere, the countries that heavily depends on fossil fuel imports will be able to use renewable energy as a manner to achieve energy independence, i.e, they will have greater energy security and more freedom to take the energy decisions that suit them, reducing its vulnerability to import fossil fuels (particularly, oil and natural gas).

However, the fast penetration of these new sources has also raised some concerns for both planners and operators that are highly studies in the literature: (i) most of these sources are non-dispatchable, i.e., their generation cannot be controlled by the system operator \cite{denholm2011grid,lund2015review,perera2017electrical}; and (ii) their energy production presents strong variability, i.e., the production can change significantly from one hour to the next \cite{halamay2011reserve,bird2013integrating,golestaneh2016very,hoeltgebaum2018generating}.

As can be seen, the VRE penetration ends up causing representative impacts on the net demand profile. In addition to the change in the profile, it is worth noting the raise of net demand ramps and their respective inclinations with the greater renewable penetration. These impacts lead to new operational challenges, which stand out: 
\begin{itemize}
    \item \textbf{Over-supply situations}: periods when the renewable generation is higher than the demand to be met (for example, in the middle of the night in regions with strong night winds or during the day in regions with a significant solar power capacity). \cite{su2017impact} in hydropower-dominated regions;
    \item \textbf{Fast upward and downward ramps}: dispatchable plants must have the ability to fast respond to the increase and decrease of intermittent renewable generation to maintain supply reliability and system stability;
    \item \textbf{Increasing thermal cycling}: possible increase in the number of startups and shutdowns of thermal plants in the system due to renewable generation intermittency;
\end{itemize}

There are several studies in the literature that adresses these challenges. For example, \cite{schaber2013managing} analyzes forms of efficiently curtail renewable generation in over-supply situations in Germany and \cite{bird2014wind} analyzes the historical operation and current practices of curtailment in the United States. Besides that, several works analyze needs for thermal flexibility due to renewable generation \cite{eser2016effect,alizadeh2016flexibility,kondziella2016flexibility}

Because of its importance, expansion planning problems are vastly discussed in the literature. There are several decomposition approaches model in the literature. In \cite{campodonico2003expansion} a Benders decomposition between investment problem and SDDP algorithm is proposed. Since this methodology is very scalable, stochastic and produces optimal solution, it has been used in several real-case studies, such as \cite{oliveira2007value,barroso2005integrated,drayton2004transmission}. The work in \cite{pina2013high} also proposes a decomposition where the master is a MILP investment problem and the slave a short-term deterministic operating model, not taking into account uncertainties. 

Since those decomposition algorithm requires convexity in the operation problem, there are some constraints such as unit commitment, that requires binary variables. In order to deal with that, there are several studies proposing the co-optimization between investment and operation, so that it could be solved with a MILP. For example, \cite{pineda2018chronological,koltsaklis2015multi,rashidaee2018linear,zhang2018mixed} proposes the co-optimization with some assumptions in order to make the MILP computational tractable. Since it may be hard to solve a huge MILP, most of the co-optimization methodology proposes to aggregate the days of the year into representative days, reducing computational time. \cite{liu2018hierarchical} showed that the clustered model (using representative days) leads to expansion planning results very similar to the unclustered model (considering all of the days of the year).


\section{Solution methodology}\label{sec:solmet}

Similar to the works in \cite{pineda2018chronological,koltsaklis2015multi,rashidaee2018linear,zhang2018mixed}, the proposed model uses the co-optimization between the investment and operation. Also, some assumptions to cluster the days of the years into representative days are used in order to reduce computational effort. Besides that, a rolling horizon scheme is implemented, in order to split the horizon into windows of a year, also in order to reduce computational effort. The next section introduces the concept of the rolling horizon,  seasons and typical days.

\subsection{Horizon Decomposition Heuristic}
\begin{figure}
    \centering
    \includegraphics[scale=0.7]{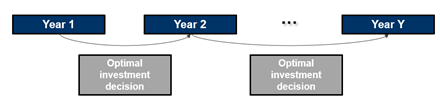}
    \caption{Horizon decomposition heuristic}
    \label{fig:hr}
\end{figure}

Since the planning horizons are long, in order to solve the expansion problem when applying co-optimization of investment and operation, the horizon is decomposed into annual sub-horizons through the forward strategy in time, that is, a problem of co-optimization of the investment and operation is solved for each year in a rolling horizon scheme. An optimum expansion plan is calculated per year, this decision is fixed, and a new optimization problem is set for the following year, considering the investment decisions taken in the previous year as fixed and completing the expansion plan, when necessary, as shown in Figure \ref{fig:hr}

\subsection{Typical days and seasons}
Since the operation is solved with hourly representation, it may result in a large and computationally intractable problem, given the size of studies that envision long-term horizons in the planning process, and since the proposed model solves a MILP that aims to minimize investment costs and the expected value of operating costs, subject to uncertainties in hydrology and generation of intermittent renewable sources. 

As a way of exemplifying this issue, taking a real energy system into account, the Table \ref{tab:comp} summarizes the size of the optimization problems for 1 month and 5 blocks versus 744 hours.

\begin{table}[]
\caption{Comparing blocks and hours resolution size of the problems}
\label{tab:comp}
\begin{tabular}{@{}|l|c|c|@{}}
\toprule
\multicolumn{1}{|c|}{\textbf{Constraints}} & \textbf{Blocks} & \textbf{Hours} \\ \midrule
Water balance constraints                  & 114             & +80,000        \\ \midrule
Load balance constraints                   & 30              & +4,000         \\ \midrule
Maximum generation \& turbining constraints & 1499            & +290,000       \\ \midrule
Maximum \& minimum volume constraints       & 228             & +165,000       \\ \midrule
Total                                      & 1461            & +520,000       \\ \bottomrule
\end{tabular}
\end{table}

As can be seen, the size of the optimization problems increases significantly. In addition to that, while evaluating real systems' expansion, it is also necessary to use multiple scenarios to incorporate the uncertainties to which the system will be exposed (hydrology, renewable generation, etc.) and, consequently, the addition of all constraints per scenario in the optimization problem. For this reason, it is necessary to create a strategy that reduces the size of the problem, but without compromising the quality of the results.

In order to reduce the computational effort required by these optimization problems, it is necessary to introduce the concepts of seasons and representative (typical) days, which in addition to enabling the solution of these problems in acceptable computational times, captures the effects of intermittent generation in the system.

\begin{figure}
    \centering
    \includegraphics[scale=0.7]{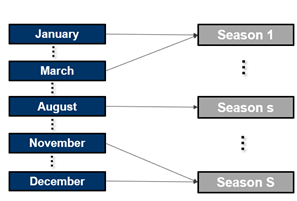}
    \caption{Mapping months into seasons}
    \label{fig:seasons}
\end{figure}

\begin{figure}
    \centering
    \includegraphics[scale=0.7]{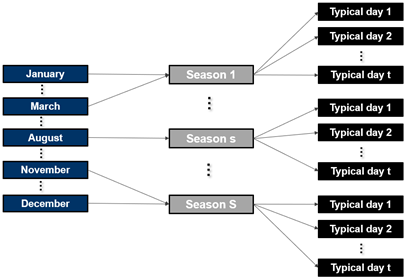}
    \caption{Mapping typical days inside each season}
    \label{fig:typicaldays}
\end{figure}

The first step of this strategy is to group the months of the year into seasons, as shown in Figure \ref{fig:seasons}. Once the seasons are defined, the representative days within each of them, here referred to as typical days, should be defined. This type of representation aims to reduce the number of days analyzed within each season, since the daily demand profiles are not usually so different, especially within the pre-defined seasons. The Figure \ref{fig:typicaldays} illustrates this grouping of real days on typical days for a set of seasons in a given year. The allocation presented in the figure was made in a generic way, with illustrative purposes.

\section{Handling uncertainties}\label{sec:uncertain}
In SDDP model, the long-term production costing decision making process (generation of each plant, interconnections between regions, circuit flows, etc.) consists in a stochastic optimization problem that seeks to balance the immediate cost and the expected value of the future cost (the expected value comes from the uncertainty about future hydrology, wind, consumption, availability of equipment, etc.). This problem is intrinsically related to storage devices that create a time-coupling between stages. Therefore, today's operating decisions, such as storage levels, may impact the mid and long-term operation, affecting thus the future operating costs. For further details, please refer to the SDDP Methodology Manual.
	
Taking the aforementioned explanation into account and given that this expansion approach performs the investment and operation co-optimization within the same problem, the operational policy is not calculated through SDDP algorithm, the proposed model does not consider the calculation of a Future Cost Function (FCF) for the system in each stage of the operation, since its calculation would require iterations of the system operation which reflects in the simulation of the operation in each stage several times until the FCF is sufficiently well approximated. It is intuitive to see that the SDDP application to calculate the FCF is the most realistic way to simulate the operation of the system, but, as it is intended to apply co-optimization, the operation of the hydro reservoirs throughout the year should be simplified.

The formulation proposed ensures that the initial storage of the reservoir of each hydroelectric plant at the beginning of each year of the study horizon will be equal to the final storage of that year. This operating strategy prevents the model from completely depleting the reservoirs present in the system during the year, optimizing its use throughout the year. The concept behind this modeling is a multi-deterministic operation, where the operation of the reservoirs is optimized for each separate scenario, without the incorporation of hydrological uncertainty into the decision-making process of the system operation in each scenario. It is plausible to explain that this simplification of the operation of large hydropower plants with large reservoirs has an optimistic bias, however, its application indicates that it is an approximation that presents satisfactory results for investment decision making and calculation of the expansion plan.

\section{Problem Formulation}\label{sec:problemform}
The expansion planning problem of an energy system is primarily formulated as a mathematical
programming problem, expressed by the formulation below. We suppose, for the sake of
simplicity, that all plants are candidate projects to the expansion problem.

\subsection{Investment Constraints}\label{sec:investmentconstraint}

\subsubsection{Precedence between projects}

\begin{align}
\label{con:precedence}
&x_{\omega} - x_{\zeta} \geq 0 &\forall \omega, \zeta\ \in P_k^{pre},\ \forall k \in R^{PRE}
\end{align}

\subsubsection{Mutually exclusive projects}

\begin{align}
\label{con:exclusivity}
&\sum\limits_{\omega \in P_{k}^{ex}}x_{\omega} \leq 1 & \forall k \in R^{ex}
\end{align}

\subsubsection{Associated projects}

\begin{align}
\label{con:associated}
&x_{\omega} - x_{\zeta} \geq 0 &\forall \omega, \zeta\ \in P_k^{as},\ \forall k \in R^{as}
\end{align}

\subsubsection{Minimum and maximum installed capacity / firm energy / firm capacity}

\begin{align}
\label{con:capacity_min}
&\sum\limits_{\zeta \in P_k^{cap}}w_\zeta^k x_{\zeta} \geq \underline{w_k} &\forall k \in R^{ctr}
\end{align}

\begin{align}
\label{con:capacity_max}
&\sum\limits_{\zeta \in P_k^{cap}}w_\zeta^k x_{\zeta} \leq \overline{w_k} &\forall k \in R^{ctr}
\end{align}

\subsection{Thermal plants constraints}
\label{sec:thermal_powerplants}

\subsubsection{Minimum and maximum energy generation}
\begin{align}
    & \underline{g_j}\gamma_{j,t,d,h,s} \leq  g_{j,t,d,h,s} \leq \overline{g_j}\gamma_{j,t,d,h,s} & \forall j,t,d,h,s
\end{align}
\subsubsection{Ramp up and ramp down generation}
\begin{align}
    & g_{j,t,d,h,s} - g_{j,t,d,h-1,s} \leq \Delta_j^{UP} & \forall j,t,d,h,s\\
    & g_{j,t,d,h-1,s} - g_{j,t,d,h,s} \leq \Delta_j^{DN} & \forall j,t,d,h,s
\end{align}
\subsubsection{Unit commitment}
\begin{align}
& st_{j,t,d,h,s} \geq \gamma_{j,t,d,h,s} - \gamma_{j,t,d,h-1,s} & \forall j,t,d,h,s\\
& \gamma_{j,t,d,h,s} \leq x_j & \forall j,t,d,h,s \label{eq:invest}\\
& \gamma_{j,t,d,h,s} \in \{0,1\} & \forall j,t,d,h,s
\end{align}	

The constraint \eqref{eq:invest} model the relation between commitment and investment decisions, preventing a thermal plant to be committed without being invested before. This constraint make continuous investment decisions to be incompatible with thermal commitment constraints (because it requires binary variables).

\subsection{Hydro plants constraints}

\subsubsection{Water storage balance}

Since the model does not consider the Future Cost to go Function (FCF), it forces water reservoir levels of all hydro plants to finish at the same level they started (initial storage = final storage), preventing the system to deplete all water in the reservoir at the end of the horizon, in order to avoid thermal operative costs. This strategy forces the model to optimize reservoir operation in order to utilize all the water inflows that arrived in the analyzed period. 

\begin{dmath}
v_{i,t+1,s} = v_{i,t,s} + a_{i,t,s} - \left( u_{i,t,s} + s_{i,t,s} \right) + \sum\limits_{m \in M_i}\left( u_{m,t,s} + s_{m,t,s} \right)\ \ \ \  \forall i,t,s 
\end{dmath}	

\begin{align}
&v_{i,T,s} = v_{i,0,s} & \forall i,s \label{con:vi=vf}
\end{align}

\subsubsection{Energy production}
The equation \eqref{eq:regulation} guarantees that the hourly energy production of the hydro plants is equal to the total energy turbined in the season. This equation assumes that the hydro plants have total regulation within season, i.e, they may freely transfer water, from an hour to another.

\begin{align}
& \sum\limits_{d,h}D_{t,d}g_{i,t,d,h,s} = \rho_i u_{i,t,s} & \forall i,t,s\label{eq:regulation}\\
&g_{i,t,d,h,s} \leq \overline{g_i}x_i & \forall i,t,d,h,s \label{con:prodmax}
\end{align}

\subsubsection{Minimum and maximum storage}
\begin{align}
&v_{i,t,s} \leq \overline{v_i}x_i & \forall i,t,s\\
&v_{i,t,s} + \delta^v_{i,t,s} = \underline{v_i}x_i & \forall i,t,s \label{con:vmax}
\end{align}
\subsubsection{Minimum and maximum turbining}
\begin{align}
&u_{i,t,s} \leq \overline{u_i}x_i & \forall i,t,s\\
&u_{i,t,s} + \delta^u_{i,t,s} = \underline{u_i}x_i & \forall i,t,s \label{con:vmax}
\end{align}
\subsubsection{Minimum total outflow}
\begin{align}
&u_{i,t,s} + s_{i,t,s} \leq \overline{q_i}x_i & \forall i,t,s\\
&u_{i,t,s} + s_{i,t,s} + \delta^q_{i,t,s} = \underline{q_i}x_i & \forall i,t,s \label{con:vmax}
\end{align}

\subsection{Renewables contraints}

Renewable plants generation decision must be lower than renewable generation scenarios.
\begin{align}
    & g_{l,t,d,h,s} \leq \psi_{l,t,d,h,s}x_l & \forall l,t,d,h,s
\end{align}

\subsection{Batteries}

\subsubsection{Energy storage balance}
Battery storage balance has hourly time steps, as in equation \eqref{con:battery_balance}. Like the hydro plants, batteries also have regulation constraints \eqref{con:battery_vivf}, where the initial energy storage is equal the final energy storage.

\begin{align}
&v_{b,t,d,h+1,s} = v_{b,t,d,h,s} + \eta^{+}_b q^{+}_{b,t,d,h,s} - q^{-}_{b,t,d,h,s} & \forall b,t,d,h,s \label{con:battery_balance}\\
&v_{b,t,d,24,s} = v_{b,t,d,0,s} & \forall b,t,d,s \label{con:battery_vivf}
\end{align}

\subsubsection{Maximum storage, charge and discharge}

\begin{align}
&v_{b,t,d,h,s} \leq \overline{V_b}x_b & \forall b,t,d,h,s \label{con:vbmax}\\
&q^{+}_{b,t,d,h,s} \leq \overline{q^{+}_b}x_b & \forall b,t,d,h,s \label{con:chmax} \\
&q^{-}_{b,t,d,h,s}s \leq \overline{q^{-}_b}x_b & \forall b,t,d,h,s \label{con:dischmax}
\end{align}

\subsection{Transmission lines constraints}

\subsubsection{Maximum flow}
The flow variables for the network representation are  $f_{k,t,d,h}^{+s}$ and $f_{k,t,d,h}^{-s}$, where these two positive variables represent the flow in both direction of each line, where $+$ means positive oriented and $-$ means negative oriented:

\begin{align}
&f_{k,t,d,h,s}^{+} \leq \overline{f_k^+}x_k & \forall k,t,d,h,s \label{con:fpmax}\\
&f_{k,t,d,h,s}^{-} \leq \overline{f_k^-}x_k & \forall k,t,d,h,s \label{con:fnmax}
\end{align}

\subsubsection{Second Kirchhoff law}
\label{sec:second_law}

The model considers two types of transmission lines: DC-Links and Circuits. Second Kirchhoff law will only be represented for circuits.

\begin{dmath}
\label{con:second_kirchoff_projects}
f_{k,t,d,h,s}^{+} - f_{k,t,d,h,s}^{-} - \gamma_k \left( \theta_{b_k^+,t,d,h,s} - \theta_{b_k^-,t,d,h,s}\right) \geq -M(1-x_k)\ \  \forall k\in K^p,t,d,h,s
\end{dmath}
\begin{dmath}
\label{con:second_kirchoff_projects}
f_{k,t,d,h,s}^{+} - f_{k,t,d,h,s}^{-} - \gamma_k \left( \theta_{b_k^+,t,d,h,s} - \theta_{b_k^-,t,d,h,s}\right) \leq M(1-x_k)\ \  \forall k\in K^p,t,d,h,s
\end{dmath}

\subsubsection{Area import/export constraints}
Area import/export constraints can limit the maximum amount of energy that enters or leave a specific electrical area.

For import constraints

\begin{align}
&\sum\limits_{k\in K^{+}_a}f_{k,t,d,h,s}^{+} + \sum\limits_{k\in K^{-}_a}f_{k,t,d,h,s}^{-} \leq \overline{Imp_a} & \forall a,t,d,h,s
\end{align}
\begin{align}
&\sum\limits_{k\in K^{+}_a}f_{k,t,d,h,s} + \sum\limits_{k\in K^{-}_a}f_{k,t,d,h,s} \geq \underline{Imp_a} & \forall a,t,d,h,s
\end{align}

For export constraints

\begin{align}
&\sum\limits_{k\in K^{-}_a}f_{k,t,d,h,s}^{+} + \sum\limits_{k\in K^{+}_a}f_{k,t,d,h,s}^{-} \leq \overline{Exp_a} & \forall a,t,d,h,s
\end{align}
\begin{align}
&\sum\limits_{k\in K^{-}_a}f_{k,t,d,h,s}^{+} + \sum\limits_{k\in K^{+}_a}f_{k,t,d,h,s}^{-} \geq \underline{Exp_a} & \forall a,t,d,h,s
\end{align}

\subsection{Generation constraint}

Generation constraint is an operative constraint which guarantees that a certain group of generators (thermal and hydro plants) always generate energy above or below a threshold.

\begin{align}
&\sum\limits_{j\in J_c^G}g_{j,t,d,h,s} + \sum\limits_{i\in I_c^G}g_{i,t,d,h,s} + \delta^g_{c,t,d,h,s} \geq \underline{g_c} & \forall c,t,d,h,s\label{con:generation_above}
\end{align}

\begin{align}
&\sum\limits_{j\in J_c^G}g_{j,t,d,h,s} + \sum\limits_{i\in I_c^G}g_{i,t,d,h,s} + \delta^g_{c,t,d,h,s} \leq \overline{g_c} & \forall c,t,d,h,s\label{con:generation_below}
\end{align}

\subsection{Reserve balance constraints}

\begin{align}
&g_{j,t,d,h,s} + r_{j,t,d,h,s} \leq \overline{g_j}\gamma_{j,t,d,h,s} & \forall j,t,d,h,s\label{con:thermal_reserve}\\
&r_{j,t,d,h,s} \leq \Delta_j^{UP} & \forall j,t,d,h,s\\
&g_{i,t,d,h,s} + r_{i,t,d,h,s} \leq \overline{g_i}x_i & \forall i,t,d,h,s \label{con:hydro_reserve}\\
&\eta_b^{-} q^{-}_{b,t,d,h,s} + r_{b,t,d,h,s} \leq \eta_b^{-} \overline{q^{-}_b}x_b & \forall b,t,d,h,s \label{con:battery_reserve1}\\
&r_{b,t,d,h,s} \leq \eta_b^{-} v_{b,t,d,h,s} & \forall b,t,d,h,s \label{con:battery_reserve2}
\end{align}

\begin{dmath}
\sum\limits_{j\in J_c^R}r_{j,t,d,h,s} + \sum\limits_{i\in I_c^R}r_{i,t,d,h,s} + \sum\limits_{b\in B_c^R}r_{b,t,d,h,s}+ \delta^R_{c,t,d,h,s} \geq R_{c,t,d,h,s} \forall c,t,d,h,s\label{con:reserve_market}
\end{dmath}

\subsection{Load balance constraints}

\begin{dmath}
\label{con:load_balance}
\sum\limits_{j\in J_n} g_{j,t,d,h,s} + \sum\limits_{i\in I_n} g_{i,t,d,h,s} + \sum\limits_{l \in L_n} g_{l,t,d,h,s} + \sum\limits_{b\in B_n}\left(\eta_b^{-} q^{-}_{b,t,d,h,s} - q^{+}_{b,t,d,h,s} \right)  + \sum\limits_{k\in K_n^{+}} \left( f_{k,t,d,h,s}^{+} - f_{k,t,d,h,s}^{-} \right) -\sum\limits_{k\in K_n^{-}}\left( f_{k,t,d,h,s}^{+} - f_{k,t,d,h,s}^{-} \right) - DE_{n,t,d,h,s} +\delta_{n,t,d,h,s} = D_{n,t,d,h,s}\ \forall n,t,d,h,s 
\end{dmath}	

\subsection{Objective function}

Let's define $\beta_{t,d,s}$ as:
\begin{equation}
    \beta_{t,d,s} = \dfrac{p_sD_{t,d}}{(1+rt)^{t-1}}
\end{equation}

Then the problem's objective function is the minimization of the following costs:
\subsubsection{Generation Cost}
\begin{dmath}
    \sum\limits_{t,d,s}\beta_{t,d,s}\left( \sum\limits_{j,h}(co_jg_{j,t,d,h,s} + cs_jst_{j,t,d,h,s}) + \sum\limits_{i,h}co_ig_{i,t,d,h,s}\right)
\end{dmath}

\subsubsection{Violation Cost}
\begin{dmath}
    \sum\limits_{i,t,s}\dfrac{p_s}{(1+rt)^{t-1}}(c_{\delta_i^v}\delta_{i,t,s}^v + c_{\delta_i^u}\delta_{i,t,s}^u + c_{\delta_i^q}\delta_{i,t,s}^q) +  \sum\limits_{t,d,s}\beta_{t,d,s}\left( \sum\limits_{c,h}(c_{\delta_c^G}\delta^G_{c,t,d,h,s} + c_{\delta_c^R}\delta^R_{c,t,d,h,s})\right)
\end{dmath}
\subsubsection{Deficit Cost}
\begin{dmath}
    \sum\limits_{t,d,s}\beta_{t,d,s}\sum\limits_{n,h}c_d\delta_{n,t,d,h,s}
\end{dmath}
\subsubsection{Elastic Demand Gain}
\begin{dmath}
    \sum\limits_{t,d,s}\beta_{t,d,s}\sum\limits_{n,h}P_n^EDE_{n,t,d,h,s}
\end{dmath}
\subsubsection{Investment Costs}
\begin{dmath}
    \sum\limits_{j\in J_x}ci_jx_j + \sum\limits_{i\in I_x}ci_ix_i + \sum\limits_{l\in L_x}ci_lx_l + \sum\limits_{b\in B_x}ci_bx_b + \sum\limits_{k\in K_x}ci_kx_k + 
\end{dmath}


\section{Conclusions}\label{sec:conclusion}
The model proposed here considers explicit operative constraints in the investment model. As a result, it can represent non-convexities in the operative constraints (such as commitment decisions). On the other hand, due to the increase of the problem's complexity, some simplifications have to be made. In this approach, we consider yearly time steps opposed to full horizon steps and representative (typical) days instead of real days within a year.

Typical days are days within a season that are considered representative of the input data. Thus, instead of representing all days of a season, the user selects a certain number of typical days to represent the season and associates these typical days with actual days. For instance, it is common to differentiate weekdays from Saturdays and Sundays, but the number of typical days and their definitions are flexible and chosen by the user.

The great advantages of this model are:
\begin{itemize}
    \item The co-optimization of investment and operating problems inside the same MILP al-lows the representation of unit commitments and other binary variables;
    \item The hourly chronological representation in the operation enables to capture the production variability of intermittent renewable sources and the generation ramps.
\end{itemize}

Besides the great advantages of this solution strategy, it's also important to remember its caveats. As explained in Section \ref{sec:uncertain}, the operative simulation is performed in a multi-deterministic way, where the operation of the reservoirs is optimized for each scenario individually, without the incorporation of hydrological uncertainty into the decision-making process of the system operation (as it is done when the SDDP methodology is applied and the FCF is calculated for each time stage). It is plausible to explain that this simplification of the operation of large hydropower plants with large reservoirs has an optimistic bias, however, its application indicates that it is an approximation that presents satisfactory results for investment decision making and calculation of the expansion plan.

Furthermore, it's also worth noting that since investment and operation problems are co-optimized in this solution strategy, then the more scenarios are contemplated in the problem, the higher computational effort will be demanded to solve the MILP. As a consequence, for large scale systems, the computational time might limit the number of scenarios that can be contemplated.

Finally, the proposed model is suitable for most real-case studies of expansion planning of renewable-dominated regions, representing hourly chronology, short-term constraints such as unit commitment and ramping, co-optimizing energy and reserves and with assumptions and approximations to make it computational tractable.

\bibliographystyle{IEEEtran}
\bibliography{Bibliography}

\begin{thebibliography}{10}
\providecommand{\url}[1]{#1}
\csname url@samestyle\endcsname
\providecommand{\newblock}{\relax}
\providecommand{\bibinfo}[2]{#2}
\providecommand{\BIBentrySTDinterwordspacing}{\spaceskip=0pt\relax}
\providecommand{\BIBentryALTinterwordstretchfactor}{4}
\providecommand{\BIBentryALTinterwordspacing}{\spaceskip=\fontdimen2\font plus
\BIBentryALTinterwordstretchfactor\fontdimen3\font minus
  \fontdimen4\font\relax}
\providecommand{\BIBforeignlanguage}[2]{{%
\expandafter\ifx\csname l@#1\endcsname\relax
\typeout{** WARNING: IEEEtran.bst: No hyphenation pattern has been}%
\typeout{** loaded for the language `#1'. Using the pattern for}%
\typeout{** the default language instead.}%
\else
\language=\csname l@#1\endcsname
\fi
#2}}
\providecommand{\BIBdecl}{\relax}
\BIBdecl

\bibitem{ondraczek2014we}
J.~Ondraczek, ``Are we there yet? improving solar pv economics and power
  planning in developing countries: The case of kenya,'' \emph{Renewable and
  Sustainable Energy Reviews}, vol.~30, pp. 604--615, 2014.

\bibitem{deichmann2011economics}
U.~Deichmann, C.~Meisner, S.~Murray, and D.~Wheeler, ``The economics of
  renewable energy expansion in rural sub-saharan africa,'' \emph{Energy
  Policy}, vol.~39, no.~1, pp. 215--227, 2011.

\bibitem{kaygusuz2012energy}
K.~Kaygusuz, ``Energy for sustainable development: A case of developing
  countries,'' \emph{Renewable and Sustainable Energy Reviews}, vol.~16, no.~2,
  pp. 1116--1126, 2012.

\bibitem{haller2012decarbonization}
M.~Haller, S.~Ludig, and N.~Bauer, ``Decarbonization scenarios for the eu and
  mena power system: Considering spatial distribution and short term dynamics
  of renewable generation,'' \emph{Energy Policy}, vol.~47, pp. 282--290, 2012.

\bibitem{oberthur2015decarbonization}
S.~Oberth{\"u}r and C.~Dupont, \emph{Decarbonization in the European Union:
  Internal policies and external strategies}.\hskip 1em plus 0.5em minus
  0.4em\relax Springer, 2015.

\bibitem{capros2014european}
P.~Capros, L.~Paroussos, P.~Fragkos, S.~Tsani, B.~Boitier, F.~Wagner, S.~Busch,
  G.~Resch, M.~Blesl, and J.~Bollen, ``European decarbonisation pathways under
  alternative technological and policy choices: A multi-model analysis,''
  \emph{Energy Strategy Reviews}, vol.~2, no. 3-4, pp. 231--245, 2014.

\bibitem{ekins2004step}
P.~Ekins, ``Step changes for decarbonising the energy system: research needs
  for renewables, energy efficiency and nuclear power,'' \emph{Energy Policy},
  vol.~32, no.~17, pp. 1891--1904, 2004.

\bibitem{optgenpsr}
\BIBentryALTinterwordspacing
PSR. (2019) Psr optgen. [Online]. Available:
  \url{https://www.psr-inc.com/softwares-en/?current=p4040}
\BIBentrySTDinterwordspacing

\bibitem{campodonico2003expansion}
N.~Campod{\'o}nico, S.~Binato, R.~Kelman, M.~Pereira, M.~Tinoco, F.~Montoya,
  M.~Zhang, and F.~Mayaki, ``Expansion planning of generation and
  interconnections under uncertainty,'' in \emph{3rd Balkans Power Conference},
  2003.

\bibitem{pereira1991multi}
M.~V. Pereira and L.~M. Pinto, ``Multi-stage stochastic optimization applied to
  energy planning,'' \emph{Mathematical programming}, vol.~52, no. 1-3, pp.
  359--375, 1991.

\bibitem{renewableenergyeurope}
M.~Tomescu, I.~Moorkens, W.~Wetzels, L.~Emele, and H.~Forster, ``Renewable
  energy in europe - approximated recent growth and knock-on effects,'' 02
  2015.

\bibitem{irenageopolitical}
I.~Global Commission on the Geopolitics~of Energy~Transformation, T.~Van~de
  Graaf, K.~Bond, and e.~al, \emph{A New World: The Geopolitics of the Energy
  Transformation}, 01 2019.

\bibitem{denholm2011grid}
P.~Denholm and M.~Hand, ``Grid flexibility and storage required to achieve very
  high penetration of variable renewable electricity,'' \emph{Energy Policy},
  vol.~39, no.~3, pp. 1817--1830, 2011.

\bibitem{lund2015review}
P.~D. Lund, J.~Lindgren, J.~Mikkola, and J.~Salpakari, ``Review of energy
  system flexibility measures to enable high levels of variable renewable
  electricity,'' \emph{Renewable and Sustainable Energy Reviews}, vol.~45, pp.
  785--807, 2015.

\bibitem{perera2017electrical}
A.~T.~D. Perera, V.~M. Nik, D.~Mauree, and J.-L. Scartezzini, ``Electrical
  hubs: An effective way to integrate non-dispatchable renewable energy sources
  with minimum impact to the grid,'' \emph{Applied Energy}, vol. 190, pp.
  232--248, 2017.

\bibitem{halamay2011reserve}
D.~A. Halamay, T.~K. Brekken, A.~Simmons, and S.~McArthur, ``Reserve
  requirement impacts of large-scale integration of wind, solar, and ocean wave
  power generation,'' \emph{IEEE Transactions on Sustainable Energy}, vol.~2,
  no.~3, pp. 321--328, 2011.

\bibitem{bird2013integrating}
L.~Bird, M.~Milligan, and D.~Lew, ``Integrating variable renewable energy:
  Challenges and solutions,'' National Renewable Energy Lab.(NREL), Golden, CO
  (United States), Tech. Rep., 2013.

\bibitem{golestaneh2016very}
F.~Golestaneh, P.~Pinson, and H.~B. Gooi, ``Very short-term nonparametric
  probabilistic forecasting of renewable energy generation - with application
  to solar energy,'' \emph{IEEE Transactions on Power Systems}, vol.~31, no.~5,
  pp. 3850--3863, 2016.

\bibitem{hoeltgebaum2018generating}
H.~Hoeltgebaum, C.~Fernandes, and A.~Street, ``Generating joint scenarios for
  renewable generation: The case for non-gaussian models with time-varying
  parameters,'' \emph{IEEE Transactions on Power Systems}, vol.~33, no.~6, pp.
  7011--7019, 2018.

\bibitem{su2017impact}
Y.~Su, J.~D. Kern, and G.~W. Characklis, ``The impact of wind power growth and
  hydrological uncertainty on financial losses from oversupply events in
  hydropower-dominated systems,'' \emph{Applied energy}, vol. 194, pp.
  172--183, 2017.

\bibitem{schaber2013managing}
K.~Schaber, F.~Steinke, and T.~Hamacher, ``Managing temporary oversupply from
  renewables efficiently: Electricity storage versus energy sector coupling in
  germany,'' in \emph{International Energy Workshop, Paris}, 2013.

\bibitem{bird2014wind}
L.~Bird, J.~Cochran, and X.~Wang, ``Wind and solar energy curtailment:
  Experience and practices in the united states,'' National Renewable Energy
  Lab.(NREL), Golden, CO (United States), Tech. Rep., 2014.

\bibitem{eser2016effect}
P.~Eser, A.~Singh, N.~Chokani, and R.~S. Abhari, ``Effect of increased
  renewables generation on operation of thermal power plants,'' \emph{Applied
  Energy}, vol. 164, pp. 723--732, 2016.

\bibitem{alizadeh2016flexibility}
M.~Alizadeh, M.~P. Moghaddam, N.~Amjady, P.~Siano, and M.~Sheikh-El-Eslami,
  ``Flexibility in future power systems with high renewable penetration: A
  review,'' \emph{Renewable and Sustainable Energy Reviews}, vol.~57, pp.
  1186--1193, 2016.

\bibitem{kondziella2016flexibility}
H.~Kondziella and T.~Bruckner, ``Flexibility requirements of renewable energy
  based electricity systems--a review of research results and methodologies,''
  \emph{Renewable and Sustainable Energy Reviews}, vol.~53, pp. 10--22, 2016.

\bibitem{oliveira2007value}
G.~C. Oliveira, S.~Binato, and M.~V. Pereira, ``Value-based transmission
  expansion planning of hydrothermal systems under uncertainty,'' \emph{IEEE
  Transactions on power systems}, vol.~22, no.~4, pp. 1429--1435, 2007.

\bibitem{barroso2005integrated}
L.~Barroso, B.~Flach, R.~Kelman, B.~Bezerra, S.~Binato, J.~Bressane, and
  M.~Pereira, ``Integrated gas-electricity adequacy planning in brazil:
  technical and economical aspects,'' in \emph{Power Engineering Society
  General Meeting, 2005. IEEE}.\hskip 1em plus 0.5em minus 0.4em\relax IEEE,
  2005, pp. 1977--1982.

\bibitem{drayton2004transmission}
G.~Drayton, M.~McCoy, M.~Pereira, E.~Cazalet, M.~Johannis, and D.~Phillips,
  ``Transmission expansion planning in the western interconnection-the planning
  process and the analytical tools that will be needed to do the job,'' in
  \emph{Power Systems Conference and Exposition, 2004. IEEE PES}.\hskip 1em
  plus 0.5em minus 0.4em\relax IEEE, 2004, pp. 1556--1561.

\bibitem{pina2013high}
A.~Pina, C.~A. Silva, and P.~Ferr{\~a}o, ``High-resolution modeling framework
  for planning electricity systems with high penetration of renewables,''
  \emph{Applied Energy}, vol. 112, pp. 215--223, 2013.

\bibitem{pineda2018chronological}
S.~Pineda and J.~M. Morales, ``Chronological time-period clustering for optimal
  capacity expansion planning with storage,'' \emph{IEEE Transactions on Power
  Systems}, 2018.

\bibitem{koltsaklis2015multi}
N.~E. Koltsaklis and M.~C. Georgiadis, ``A multi-period, multi-regional
  generation expansion planning model incorporating unit commitment
  constraints,'' \emph{Applied energy}, vol. 158, pp. 310--331, 2015.

\bibitem{rashidaee2018linear}
S.~A. Rashidaee, T.~Amraee, and M.~Fotuhi-Firuzabad, ``A linear model for
  dynamic generation expansion planning considering loss of load probability,''
  \emph{IEEE Transactions on Power Systems}, vol.~33, no.~6, pp. 6924--6934,
  2018.

\bibitem{zhang2018mixed}
Y.~Zhang, Y.~Hu, J.~Ma, and Z.~Bie, ``A mixed-integer linear programming
  approach to security-constrained co-optimization expansion planning of
  natural gas and electricity transmission systems,'' \emph{IEEE Transactions
  on Power Systems}, 2018.

\bibitem{liu2018hierarchical}
Y.~Liu, R.~Sioshansi, and A.~J. Conejo, ``Hierarchical clustering to find
  representative operating periods for capacity-expansion modeling,''
  \emph{IEEE Transactions on Power Systems}, vol.~33, no.~3, pp. 3029--3039,
  2018.

\end{thebibliography}

\begin{IEEEbiography}[{\includegraphics[width=1in,height=1.25in,clip,keepaspectratio]{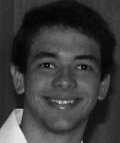}}]%
{Alessandro Soares}
Has a BSc in Electrical Engineering and in Control Engineering from PUC-Rio. Is currently doing a MSc in Optimization/Operations Research at PUC-RJ. Joined PSR in 2017 and has been working on the development and support of the expansion planning model (OPTGEN) and with the time series analysis model Time Series Lab (TSL). Before joining PSR, worked with non-linear time series models for synthetic inflow scenarios generation and modelling probabilistic scenarios for renewable plants. He has also worked with new methodologies for the analysis of the impact of climate changes in the inflow series.
\end{IEEEbiography}

\begin{IEEEbiography}[{\includegraphics[width=1in,height=1.25in,clip,keepaspectratio]{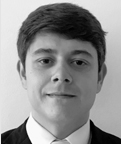}}]%
{Ricardo Perez}
Has a BSc in Electrical Engineering with emphasis in Power Systems from Itajubá Federal University (UNIFEI). During college he worked in two R\&D; projects in the Power Quality Study Group. Through an exchange program with the Technische Universitat Dresden in Germany, he also developed a research in the same field at this university. In addition to the applied research, his experience in the electricity business includes internships in Brazil and in Germany, respectively at the Generation Planning Department of CPFL Geracao and at DIgSILENT GmbH (where he carried out studies regarding the connection of Wind Farms to the German grid). Mr. Perez joined PSR in December 2009 and has been a member of the transmission studies group ever since.
\end{IEEEbiography}

\begin{IEEEbiography}[{\includegraphics[width=1in,height=1.25in,clip,keepaspectratio]{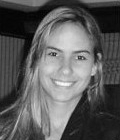}}]%
{Fernanda Thome}
Has a BSc in EE from UFRJ and MSc in OR from COPPE/UFRJ, where she is currently working towards a DSc in OR. She joined PSR in 2000 where she initially took part in the development of transmission expansion databases using geographical information systems, development of computational tools applied to the optimization models for operation and expansion planning of electric power systems, energy pricing model and viability studies on transmission systems expansion. Recently she has been participating in studies and development of hydrothermal operation planning model and expansion planning models for generation systems, transmission and networks and natural gas networks.
\end{IEEEbiography}

\end{document}